\documentclass[a4paper, 10pt]{article}
\usepackage{mathrsfs}

\usepackage{epsfig}
\usepackage{amsmath}
\usepackage{amssymb}
\usepackage{latexsym}
\usepackage{amsfonts}
\usepackage{amsthm}
\usepackage[numbers,sort&compress]{natbib}
\usepackage{graphicx}
\ifx\pdfoutput\undefined  \else
 \fi
\usepackage[centerlast]{caption2}
\usepackage{color}

\newtheorem{thm}{Theorem}[section] 

\newtheorem{obs}[thm]{Observation}

\newtheorem{lemma}[thm]{Lemma}

\newtheorem{example}[thm]{Example}
\newtheorem{conj}[thm]{Conjecture}

\newtheorem{defi}[thm]{Definition}

\def\qed{\nopagebreak\hfill{\rule{4pt}{7pt}}\medbreak}

\allowdisplaybreaks

\begin{document}

\title{{Extremal incomplete sets in \\ finite abelian groups}
 \footnote{The research is
supported by NSFC (11001035)}
\author{
Dan Guo$^{a}$,\ Yongke Qu$^{a,}$\footnote {corresponding email:
quyongke@mail.nankai.edu.cn}, \ Guoqing Wang$^{a}$,\ Qinghong Wang$^{b}$ \\
\\
$^{a}$Center for Combinatorics, LPMC-TJKLC\\
Nankai University\\
Tianjin 300071, P. R. China\\
\\
$^{b}$College of Science\\
Tianjin University of Technology\\
Tianjin 300384, P. R. China }}


\date{}
\maketitle

\begin{abstract}
Let $G$ be a finite abelian group. The critical number ${\rm cr}(G)$
of $G$ is the least positive integer $\ell$ such that every subset
$A\subseteq G\setminus\{0\}$ of cardinality at least $\ell$ spans
$G$, i.e., every element of $G$ can be written as a nonempty sum of
distinct elements of $A$. The exact values of the critical number
have been completely determined recently for all finite abelian
groups. The structure of these sets of cardinality ${\rm cr}(G)-1$
which fail to span $G$ has also been characterized except for the
case that $|G|$  is an even number and the case that $|G|=pq$ with
$p,q$ are primes. In this paper, we characterize these extremal
subsets for $|G|\geq 36$ is an even number, or $|G|=pq$ with $p,q$
are primes and $q\geq 2p+3$.
\end{abstract}

\section{Introduction}

Let $G$ be a finite abelian group, and let $A$ be a subset of $G$.
Let $\langle A\rangle$ denote the subgroup of $G$ generated by $A$.
Let $$\sum(A)=\{\sum\limits_{g\in B} g: \emptyset\neq B\subseteq
A\}.$$ The nonempty set $A$ is said to be {\sl complete} if
$\sum(A)=\langle A\rangle$. We say $A$ spans $G$ if $\sum(A)=G$,
equivalently, if $\langle A\rangle=G$ and $A$ is complete. The
critical number ${\rm cr}(G)$ of
 $G$ is the least positive integer $\ell$ such that every
subset $A\subseteq G\setminus\{0\}$ of cardinality at least $\ell$
spans $G$. The critical number was first studied by P. Erd\H{o}s and
H. Heilbronn \cite{Erdos} for cyclic groups of prime order in 1964.
Since then, due to contributions by H.B. Mann, J.E. Olson, G.T.
Dierrich, Y.F. Wou, J.A. Dias da Silva,  Y.O. Hamidoune, A.S.
Llad\'{o}, O. Serra, M. Freeze, W.D. Gao, and A. Geroldinger,  et
al., the critical numbers of all finite abelian groups have been
completely determined. Also, there has been some work generalizing
the critical number to non-commutative groups, one can refer to
\cite{GaoNoncom,WangQu}. The values of critical number of finite
abelian groups are summarized as follows.

{\bf Theorem A.}  (\cite{h}, \cite{Did75}, \cite{Fre-Gao-Ger},
\cite{Gao-Ham}, \cite{Gri})\ {\sl Let $G$ be a finite abelian group
of order $|G|\geq 3$, and let $p$ denote the smallest prime divisor
of $|G|$.

1. If $|G|=p$, then ${\rm cr}(G)=\lfloor2\sqrt{p-2}\rfloor$.

2. In each of the following cases we have ${\rm
cr}(G)=\frac{|G|}{p}+p-1$;

$\centerdot$ G is isomorphic to one of the following groups:
$Z_2\oplus Z_2$, $Z_3\oplus Z_3$, $Z_4$, $Z_6$, $Z_2\oplus Z_4$,
$Z_8$.

$\centerdot$ $\frac{|G|}{p}$ is an odd prime with
$2<p<\frac{|G|}{p}\leq p+\lfloor2\sqrt{p-2}\rfloor+1$.

3. In all other cases we have ${\rm cr}(G)=\frac{|G|}{p}+p-2$.}

So, a natural question is

{\bf What is the structure of the extremal subsets which fail
to span $G$?}\\

\bigskip

With respect to this question, Nguyen, Szemed\'{e}di, and Vu
\cite{Vu} characterized the set of $Z_p$ of cardinality at least
$\sqrt{2p}$ which fails to span the group $Z_p$. Recently, Vu
\cite{Vu2} also showed that for general finite abelian group $G$
meeting certain conditions, if
 $A$ is a comparatively large subset of $G$ and $A$ fails to span $G$,
then $A$ contains a complete subset. Before then, Gao, Hamidoune,
Llad\'{o} and Serra \cite{G-H-L-S} obtained the following result:

{\bf Theorem B.} \  {\sl Let $G$ be a finite abelian group of odd
order. Let $p$ be the smallest prime divisor of $|G|$. Assume
$\frac{|G|}{p}$ is composite and
$$\begin{array}{llll} & \frac{|G|}{p}\geq\left \{\begin{array}{llll}
               62, & \mbox{ if } \ \ p=3;\\
               7p+3, & \mbox{ if } \ \ p\geq 5.\\
              \end{array}
           \right. \\
\end{array}$$
Let $A$ be a subset of $G\setminus \{0\}$ of cardinality ${\rm
cr}(G)-1$ such that $\sum(A)\neq G$. Then there exists a subgroup
$H$ of order $\frac{|G|}{p}$ and an element $g\in G\setminus H$ such
that $H\setminus \{0\}\subseteq A$ and $A\subseteq H\cup (g+H)\cup
(-g+H)$.}

In general, the structure of the extremal set remains unknown only
for the following two types of group $G$:

1. $|G|$ is even.

2. $|G|$ is a product of two odd prime numbers.

In this paper, we characterize the structure of the extremal set for
the group $G$ when $|G|$ is an even number with $|G|\geq 36$, or
$|G|$ is a product of two odd prime numbers $p,q$ with $q\geq 2p+3$.
Our main result is as follows.

\begin{thm}\label{Theorem main}\ Let $G$ be a finite abelian group, and let $p$ be the smallest prime
divisor of $|G|$. Assume that $p=2$ and $|G|\geq 36$, or
$\frac{|G|}{p}$ is a prime number with $\frac{|G|}{p}\geq 2p+3$. Let
$A\subseteq G\setminus \{0\}$ be a subset of cardinality ${\rm
cr}(G)-1$ such that $\sum(A)\neq G$. Then there exists a subgroup
$H$ of $G$ of cardinality $\frac{|G|}{p}$ such that
\begin{itemize}
\item[(i)] \ \ if $p=2$, then $A=H\setminus \{0\}$;
\item[(ii)] \ \ if $p\geq3$, then $A\subseteq H\cup (g+H)\cup
(-g+H)$ and $H\setminus \{0\}\subseteq A$, where $g\in G\setminus
H$.
\end{itemize}
\end{thm}

The rest of this paper is organized as follows. In Section 2 we
introduce some technical notations and tools. The proof of Theorem
\ref{Theorem main} is presented in Section 3. The Final Section 4
contains some concluding remarks, together with two conjectures on
the structure of the extremal set for the group $G=Z_{pq}$ with
$p<q< 2p+3$.

\section{Notations and tools}

Let $G$ be a finite abelian group, and let $A$ and $B$ be nonempty
subsets of $G$. The sumset $A+B$ is the set of all elements of $G$
that can be written in the form $a+b$, where $a\in A$ and $b\in B$.
We call $A$ an {\sl arithmetic progression with difference $d\in G$}
if there is some $a\in G$ such that $A=\{a+id: i\in [0,|A|-1]\}$.
Let $T=a_1a_2\cdot\ldots\cdot a_{\ell}$ be a sequence of elements in
$G$. Define $$\sum(T)=\{\sum\limits_{i\in I}a_i: \emptyset\neq
I\subseteq [1,\ell]\}$$ and
$$\sum\nolimits^{\circ}(T)=\sum(T)\cup \{0\}.$$
For notational convenience, we let $\sum(T)=\{0\}$ if $T$ is an
empty sequence. For any integer $h\in [1,\ell]$,  let
$$\sum\nolimits_h(T)=\{\sum\limits_{i\in I}a_i: I\subseteq [1,\ell], |I|=h\}.$$
We adopt the convention that $\sum\nolimits_0(T)=\{0\}$. Let ${\rm
supp}(T)$ be the set consisting of all distinct elements in $T$. In
this paper, we shall view a set $A$ to be a squarefree sequence,
i.e., ${\rm supp}(A)=A$. Then all the notations that are valid for
sequences automatically apply to sets too.

We present below some tools:

\begin{lemma}[\cite{Geroldinger1}, Lemma 5.2.9]\label{Lemma folk}\  Let $A_1$ and $A_2$ be nonempty subsets of a finite abelian group $G$.
If $|A_1|+|A_2|\geq |G|+1$ then $A_1+A_2=G$.
\end{lemma}

\begin{lemma}[\cite{Hamidoune99}]\label{Lemma Hamidoune}\  Let $A$ be a subset of a finite abelian group $G$ such that
$0\notin A$ and $|A|\geq 14$. Then one of the following conditions
holds.
\begin{enumerate}
\renewcommand{\labelenumi}{(\roman{enumi}).}
\item $|\sum\nolimits^{\circ}(A)|\geq \min(|G|-3,3|A|-3)$;
\item There is a subgroup $H\neq G$ such that $|A\cap H|\geq |A|-1$.
\end{enumerate}
\end{lemma}

\begin{lemma}[\cite{Nat96}, Theorem 2.3]\label{Lemma Cauchy} \ Let $p$ be a prime number,
and let $A_1,A_2,\ldots,A_h$ be nonempty subsets of $Z_p$. Then
$$|A_1+A_2+\cdots+A_h| \geq \min(p, \sum\limits_{i=1}^h |A_i|-h+1).$$
\end{lemma}

\begin{lemma}[\cite{Did75}]\label{Lemma Diderrich}  \ Let $p$ be a prime number, and let $A_1,A_2,\ldots,A_h$ be nonempty subsets of
$Z_p$, apart from one possible exception, are arithmetic
progressions with pairwise distinct nonzero differences. Then
$$|A_1+A_2+\cdots +A_h|\geq\min(p,\sum\limits_{i=1}^h|
A_i|-1).$$
\end{lemma}

\noindent{\bf Remark.} Note that an arithmetic progression
$A=\{a,a+d,\ldots,a+(|A|-1)d\}$ of difference $d$ can  also be
viewed as an arithmetic progression of difference $-d$. Let
$A_1,A_2,\ldots,A_h$ be $h$ arithmetic progressions with
$A_i=\{a_i,a_i+d_i,\ldots,a_i+(|A_i|-1)d_i\}$. If one can find some
$h-$tuple $(\theta_1,\theta_2,\ldots,\theta_h)\in
\{d_1,-d_1\}\times\{d_2,-d_2\}\times \cdots \times \{d_h,-d_h\}$
such that $\theta_1, \theta_2, \ldots, \\ \theta_h$ are pairwise distinct,
then $A_1,A_2,\ldots,A_h$ would be regarded as progressions with
pairwise distinct differences.

\begin{lemma} [\cite{Nat96}, Theorem 2.7]\label{Lemma Vosper}\  Let $p$ be an odd prime number, and let $B_1$ and $B_2$ be
nonempty subsets of the group $Z_p$ with $2\leq |B_i|\leq p-2$ for
$i=1,2$. If $|B_1+B_2|<\min(p,|B_1|+|B_2|)$, then $B_1$ and $B_2$
are arithmetic progressions with differences $d_1$ and $d_2$,
respectively, such that $d_1\in \{d_2,-d_2\}$.
\end{lemma}

\begin{lemma}\label{lemma three facts} \ Let $p$ be a prime number,
and let $A$ be a nonempty subset of $Z_p$. Then
\begin{enumerate}
\renewcommand{\labelenumi}{(\roman{enumi}).}
\item $|\sum\nolimits_h(A)|\geq \min(p,h|A|-h^2+1)$ for all $1\leq h\leq |A|$
(see \cite{h}, \cite[Theorem 3.4]{Nat96});
\item If $|A|=\lfloor\sqrt{4p-7}\rfloor$ and $h=\lfloor|A|/2\rfloor$,
then $|\sum\nolimits_h(A)|=p$ (see \cite{h});
\item If $|A|\geq \lfloor2\sqrt{p-2}\rfloor$ then $|\sum(A)|=p$ (see \cite{Fre-Gao-Ger}).
\end{enumerate}
\end{lemma}

The following lemma is a corollary of Theorem 1.3 in
\cite{GaoPengWang}.

\begin{lemma}\label{Lemma from behaving}\  Let $A$ be a nonempty, finite subset of an abelian group with $0\notin
A$. Then $$|\sum(A)|\geq \min(|\langle A\rangle|, 2|A|-1).$$
\end{lemma}

By Lemma \ref{Lemma from behaving}, we immediately have the
following

\begin{lemma}\label{Lemma case as ell0}\  Let $p$ be a prime number, and let $A$ be a subset of
$Z_p\setminus \{0\}$ with $|A|=\ell$. Then
$$|\sum\nolimits^{\circ}(A)|\geq \min (p,2\ell-1+\epsilon(\ell)),$$
where
$$\begin{array}{llll} & \epsilon(\ell)=\left \{\begin{array}{llll}
               2, & \mbox{ if } \ \ \ell=0;\\
               1, & \mbox{ if } \ \ \ell=1;\\
               0, & \mbox{ if } \ \ \ell\geq 2.\\
              \end{array}
           \right. \\
\end{array}$$
\end{lemma}

\begin{lemma}\label{Lemma fact on sequence sums}\
Let $p$ be a prime number, and let $T$ be a sequence of elements in
$Z_p\setminus \{0\}$ and of length $|T|\geq 2$. Then
$|\sum\nolimits^{\circ}(T)|\geq \min(p,|T|+1)$, and moreover,
equality holds only for one of the following two conditions.

\begin{enumerate}
\renewcommand{\labelenumi}{(\roman{enumi}).}
\item $|T|\geq p-1$;
\item There exists some $g\in Z_p\setminus \{0\}$ such that ${\rm
supp}(T)\subseteq \{g,-g\}$.
\end{enumerate}
\end{lemma}

\begin{proof} Let $T=a_1\cdot \ldots\cdot a_{\ell}$.
Let $A_i=\{0,a_i\}$ for $i\in [1,\ell]$. By Lemma \ref{Lemma
Cauchy}, we have
\begin{align*}
|\sum\nolimits^{\circ}(T)|&= |\sum\limits_{i=1}^{\ell}A_i| \geq
\min(p, \sum\limits_{i=1}^{\ell}|A_i|-{\ell}+1) =\min(p,\ell+1).
\end{align*}
Now assume that neither (i) nor (ii) holds, i.e.,
$$\ell\leq p-2$$ and there exist two elements, say $$a_{\ell-1}\neq
a_{\ell}$$ and $$a_{\ell-1}+a_{\ell}\neq 0.$$ Let
$A_0=\{0,a_{\ell-1},a_{\ell},a_{\ell-1}+a_{\ell}\}$. It follows from
Lemma \ref{Lemma Cauchy} that
\begin{align*}
|\sum\nolimits^{\circ}(T)|&= |A_0+A_1+\cdots +A_{\ell-2}|\geq
\min(p,\sum\limits_{i=0}^{\ell-2}|A_i|-(\ell-2)) =\ell+2.
\end{align*} Then the lemma follows.
\end{proof}

\section{Proof of Theorem \ref{Theorem main}}

We begin this section with the following observation.

\begin{obs}\label{Observation cap one subgroup}\  Let $G$ be a finite abelian group,
and let $A$ be a subset of $G\setminus \{0\}$ of cardinality ${\rm
cr}(G)-1$ such that $\sum(A)\neq G$. If $A^{'}$ is a complete subset
of $A$, then $A\cap \langle A^{'}\rangle=\langle
A^{'}\rangle\setminus \{0\}$.
\end{obs}

We shall prove Theorem \ref{Theorem main} by two cases according to
$|G|$ is an even number or $|G|$ is a product of two distinct prime
numbers.
\\

\bigskip

{\sl Proof of Theorem \ref{Theorem main} for the case that
$|G|\equiv 0\pmod 2$ with $|G|\geq 36$.}

By Theorem A, we have
\begin{equation}\label{equation |A| geq 17}
|A|={\rm cr}(G)-1=\frac{|G|}{2}-1\geq 17.
\end{equation}
 Take a subset $A_1$ of $A$
with
\begin{equation}\label{equation |A1|=3 in even number}
|A_1|=3.
\end{equation} Let $A_2=A\setminus A_1$. Obviously,
$|\sum(A_1)|\geq 4.$ Since
$\sum(A_1)+\sum\nolimits^{\circ}(A_2)\subseteq \sum(A)\neq G$, it
follows from Lemma \ref{Lemma folk} that
$|\sum\nolimits^{\circ}(A_2)|<|G|-3$. By Lemma \ref{Lemma
Hamidoune}, we conclude that there is a subgroup $K\neq G$ such that
$|A_2\cap K|\geq |A_2|-1$. It follows from \eqref{equation |A| geq
17} and \eqref{equation |A1|=3 in even number} that $|A_2\cap K|\geq
|A_2|-1=\frac{|G|}{2}-4>\frac{|G|}{3}-1$, and so
$$|K|=\frac{|G|}{2}.$$
By \eqref{equation |A| geq 17}, we have $|A_2\cap K|\geq |K|-4\geq
\frac{|K|+1}{2}.$ It follows from Lemma \ref{Lemma from behaving}
that $\sum(A_2\cap K)=\langle A_2\cap K\rangle=K$. Then the
conclusion follows from Observation \ref{Observation cap one
subgroup}. \qed

Therefore, it remains to prove Theorem \ref{Theorem main} for the
case that $$G=Z_{pq}$$ where $p,q$ are odd prime number such that
$$q\geq  2p+3.$$
Before proceeding with our arguments, we need to formulate some more
technical notations and definitions which will be used in the rest
part of this paper.

Let $H$ be the subgroup of $G$ of order $q$, and let $\varphi$ be
the canonical epimorphism of $G$ onto the quotient group $G\diagup
H$. Then $\varphi(A)$ is a sequence of elements in $G\diagup H$  of
length $p+q-3$. Denote $$k=|{\rm supp}(\varphi(A))\setminus
\{0\}|.$$ Fix $k$ elements $a_1,\ldots,a_k\in A\setminus H$ such
that $\varphi(a_1),\ldots, \varphi(a_k)$ are pairwise distinct. Let
$$A\setminus H=\bigcup\limits_{i=1}^k A_i,$$ where $A_i+H=a_i+H$ for all $i\in [1,k]$, and let
$$A_0=A\cap H$$ ($A_0$ is perhaps an empty set). Let
$$\widetilde{A_i}=A_i-a_i\ \ \ \mbox{ for }\ i\in [1,k].$$
Note that $\widetilde{A_1},\ldots, \widetilde{A_k}$ are subsets of
$H$. Denote $${\ell}_i=|A_i|\ \ \ \mbox{ for }\ i\in [0,k].$$  We
shall always admit
\begin{equation}\label{equation ell is decreasing}
\ell_1\geq \ell_2\geq \cdots \geq \ell_k.
\end{equation}
Let
\begin{align*}
&R_1=\{ i\in [1,k]: \ell_i=1\};\\
&R_2=\{i\in [1,k]: \ell_i=2\};\\
&R_3=\{i\in [1,k]: \ell_i=3\}; \\
&R_4=\{i\in [1,k]: \ell_i=4\}; \\
&R_5=\{i\in [1,k]: \ell_i\geq 5\}; \\
\end{align*}
and let $$r_i=|R_i| \ \ \  \mbox{ for }   i\in [1,5].$$ For
convenience, let
$$m_t=k-\sum\limits_{i=1}^{t-1}r_i\ \ \  \mbox{ for } t\in [1,5].$$
Notice that $$[1,m_t]=\{i\in [1,k]:\ell_i\geq t\}\ \ \  \mbox{ for }
t\in [1,5],$$ and that
$$[m_{u+1}+1,m_u]=\{i\in [1,k]:\ell_i=u\}\ \ \  \mbox{ for } u\in
[1,4].$$

\begin{defi} \ For any element $\bar{g}\in G\diagup H$, we say that $\bar{g}$ has a representation
with coefficients $f_1,\ldots,f_k$ provided that
\begin{equation}\label{equation f1,f2,...}
\bar{g}=\sum\nolimits_{i=1}^k f_i \ \varphi(a_i),
\end{equation}
where $f_i\in [0,\ell_i]$ and $f_1+\cdots +f_k>0.$
\end{defi}

\begin{defi} \ Let
\begin{align*}
X_i^{\Delta}&=\{2\varphi(a_i),3\varphi(a_i),\ldots,
\mu_i\varphi(a_i)\} \ \ \ \mbox{ if }\ \  i\in [1,m_4];\\
X_i&= \{\varphi(a_i),2\varphi(a_i), \ldots, \lambda_i\varphi(a_i)\}
\ \ \ \ \mbox{ if }\ \  i\in [1,m_2];\\
X_i&=\{0,\varphi(a_i)\} \ \  \ \ \ \ \ \ \ \ \ \ \ \ \ \ \ \ \ \ \ \ \ \ \ \mbox{ if }\ \  i\in [m_2+1,k];\\
Y&=\varphi(\sum\limits_{i=1}^{r_2}
a_{m_3+i})-\{0,\varphi(a_{m_3+1}),
\varphi(a_{m_3+2}),\ldots,\varphi(a_{m_3+r_2})\},
\end{align*}
where  $\mu_i=\min(\ell_i-2,p+1)$ and $\lambda_i=\min(\ell_i-1,p)$.
\end{defi}

\begin{defi}\  For $n\in [0,m_4]$,  let
$$\alpha_n=\min(p,|\sum\limits_{i=1}^n
X_i^{\Delta}+\sum\limits_{i=n+1}^{m_3} X_i+\sum\limits_{i=m_2+1}^k
X_i+Y|).$$
\end{defi}

For any integer $m\geq 0$, let
$$\begin{array}{llll} & \delta(m)=\left \{\begin{array}{llll}
               0, & \mbox{ if } \ \ m=0;\\
               1, & \mbox{ if } \ \ m>0.\\
              \end{array}
           \right. \\
\end{array}$$

\bigskip

To make the remainder of the proof clear, we propose the general
idea as follows. If $|A_0|$ is large, the conclusion of Theorem
\ref{Theorem main} is easy to prove. Assume $|A_0|$ is small. We
shall derive a contradiction by the following process. We first show
that there exists some $n\in [0,m_4]$ such that $\alpha_n=p$, i.e.,
every element $\bar{g}\in G\diagup H$ has a fixed representation
with coefficients $f_1=f_1(n,\bar{g}),\ldots, f_k=f_k(n,\bar{g})$,
where
\begin{equation}\label{equation f1,...,fn}
f_i\in [2,\ell_i-2] \ \ \ \ \mbox{ for } i\in [1,n],
\end{equation}
and that
 \begin{equation}\label{equation fn+1...} f_j\in
[1,\ell_j-1] \ \ \ \ \mbox{ for } j\in [n+1,m_2]
\end{equation}
with at most one exception $u\in R_2$ such that $f_u=0$,\\
 and that
\begin{equation}\label{equation fi i in m2+1..} f_w\in
[0,1] \ \ \ \ \mbox{ for } w\in [m_2+1,k].
\end{equation}
Furthermore, we show that
$$|\sum\nolimits^{\circ}(A_0)+\sum\nolimits_{f_1}(\widetilde{A_1})+
\cdots+\sum\nolimits_{f_k}(\widetilde{A_k})|=q$$ for every
$\bar{g}\in G\diagup H$, where $f_1=f_1(n,\bar{g}),\ldots,
f_k=f_k(n,\bar{g})$ are given as \eqref{equation f1,...,fn},
\eqref{equation fn+1...} and \eqref{equation fi i in m2+1..}. This
would implies $\sum(A)=G$, which is a contradiction.

Based on the above, we shall require the following two lemmas.

\begin{lemma}\label{Lemma alpha n}\  For $n\in [0,m_4]$, if $$q-3-\ell_0-3n+\delta(r_2)\geq
p$$ then $\alpha_n=p$.
\end{lemma}

\begin{proof} We may assume without loss of generality that
$\ell_i\leq p+2$ for $i\in [1,n]$, and that $\ell_j\leq p$ for $j\in
[n+1,m_4]$. This implies that
\begin{equation}\label{equation Xidelta=ell i-3}
|X_i^{\Delta}|=\ell_i-3 \ \ \ \mbox{ for }\ i\in [1,n],
\end{equation}
and
\begin{equation}\label{equation Xj=ell j-1}
|X_j|=\ell_j-1 \ \ \ \mbox{ for }\ j\in [n+1,m_3].
\end{equation}
Denote $t=|[1,n]\cap R_4|$. That is,
\begin{equation}\label{equation ell n-t+1=...=4}
\ell_{n-t+1}=\ell_{n-t+2}=\cdots=\ell_n=4
\end{equation}
 and
$$\ell_i\geq 5\ \ \ \mbox{ for }\ i\in [1,n-t].$$
Thus,
\begin{equation}\label{equation X_i delta =1}
|X_i^{\Delta}|=1 \ \ \ \mbox{ for }\ i\in [n-t+1,n].
\end{equation}\\
It follows from \eqref{equation Xidelta=ell i-3}, \eqref{equation
Xj=ell j-1}, \eqref{equation ell n-t+1=...=4}, \eqref{equation X_i
delta =1} and Lemma \ref{Lemma Diderrich} that
\begin{eqnarray*}
\alpha_n&=&\min(p,|\sum\limits_{i=1}^{n-t}
X_i^{\Delta}+\sum\limits_{i=n-t+1}^n X_i^{\Delta}+
\sum\limits_{i=n+1}^{m_3} X_i+
\sum\limits_{i=m_2+1}^{k}X_i+Y|)\\
&=&\min(p,|\sum\limits_{i=1}^{n-t} X_i^{\Delta}+
\sum\limits_{i=n+1}^{m_3} X_i+
\sum\limits_{i=m_2+1}^{k}X_i+Y|)\\
&=&\min(p,\sum\limits_{i=1}^{n-t} |X_i^{\Delta}|+
\sum\limits_{i=n+1}^{m_3} |X_i|+
\sum\limits_{i=m_2+1}^{k}|X_i|+(r_2+\delta(r_2))-1)\\
&=&\min(p,\sum\limits_{i=1}^{n-t} (\ell_i-3)+
\sum\limits_{i=n+1}^{m_3} (\ell_i-1)+
\sum\limits_{i=m_2+1}^{k}(\ell_i+1)\\
& &\ \ \ \ \ \ \ \ \ +(\sum\limits_{i=m_3+1}^{m_2}(\ell_i-1)+\delta(r_2))-1)\\
&=&\min(p,\sum\limits_{i=1}^k \ell_i-
\sum\limits_{i=n-t+1}^{n} \ell_i-3(n-t)-(m_3-n)\\
& &\ \ \ \ \ \ \ \ \ +(k-m_2)-(m_2-m_3)+\delta(r_2)-1)\\
&=&\min(p,(p+q-3-\ell_0)-
4t-3(n-t)+n+k\\
& &\ \ \ \ \ \ \ \ \ -2m_2+\delta(r_2)-1)\\
&=&\min(p,p+q-3-\ell_0-
t-2n+k-2(k-r_1)+\delta(r_2)-1)\\
&\geq&\min(p,p+q-3-\ell_0-
n-2n+k-2(k-r_1)+\delta(r_2)-1)\\
   &\geq&\min(p,p+q-3-\ell_0-3n-k+\delta(r_2)-1)\\
   &\geq&\min(p,p+q-3-\ell_0-3n-(p-1)+\delta(r_2)-1)\\
&=&p.\\
\end{eqnarray*}
Then the lemma follows.
\end{proof}

\begin{lemma}\label{Lemma alpha and beta} \ Assume that there exists
 some $n\in [0,m_4]$ such that $\alpha_n=p$ and $\sum\limits_{i=0}^n\ell_i+
 \epsilon(\ell_0)-3n-\delta(r_2)\geq 3$. For every element $\bar{g}
 \in G\diagup H$, $$|\sum\nolimits^{\circ}(A_0)+\sum\limits_{i=1}^k
\sum\nolimits_{f_i}(\widetilde{A_i})|=q$$  where
 $f_1=f_1(n,\bar{g}),\ldots,f_k=f_k(n,\bar{g})$.
\end{lemma}

\begin{proof}

By \eqref{equation f1,...,fn}, \eqref{equation fn+1...},
\eqref{equation fi i in m2+1..}, Lemma \ref{Lemma Cauchy}, Lemma
\ref{lemma three facts} (i) and Lemma \ref{Lemma case as ell0}, we
have that
\begin{eqnarray*}
& &|\sum\nolimits^{\circ}(A_0)+\sum\limits_{i=1}^k \sum\nolimits_{f_i}(\widetilde{A_i})|\\
&=&|\sum\nolimits^{\circ}(A_0)+\sum\limits_{i=1}^{m_2}\sum\nolimits_{f_i}(\widetilde{A_i})|\\
&\geq&\min(q,|\sum\nolimits^{\circ}(A_0)|+\sum\limits_{i=1}^n|\sum\nolimits_{f_i}(\widetilde{A_i})|
+\sum\limits_{i=n+1}^{m_3}|\sum\nolimits_{f_i}(\widetilde{A_i})|\\
& &\ \ \ \ \ \ \ \ \ +\sum\limits_{i=m_3+1}^{m_2}|\sum\nolimits_{f_i}(\widetilde{A_i})|-m_2)\\
&\geq&\min(q,(2\ell_0-1+\epsilon(\ell_0))+\sum\limits_{i=1}^n(2\ell_i-3)+\sum\limits_{i=n+1}^{m_3}\ell_i\\
& &\ \ \ \ \ \ \ \ \ +(\sum\limits_{i=m_3+1}^{m_2}\ell_i-\delta(r_2))-(k-r_1))\\
&\geq&\min(q,(2\ell_0-1+\epsilon(\ell_0))+\sum\limits_{i=1}^n(2\ell_i-3)+\sum\limits_{i=n+1}^{m_3}\ell_i\\
& &\ \ \ \ \ \ \ \ \ +(\sum\limits_{i=m_3+1}^{m_2}\ell_i-\delta(r_2))-(k-\sum\limits_{i=m_2+1}^{k}\ell_i))\\
&=&\min(q,\sum\limits_{i=0}^k\ell_i+\sum\limits_{i=0}^n \ell_i+\epsilon(\ell_0)-3n-\delta(r_2)-(k+1))\\
&\geq&\min(q,p+q-3+\sum\limits_{i=0}^n \ell_i+\epsilon(\ell_0)-3n-\delta(r_2)-p)\\
&=&q.
\end{eqnarray*}

Then the lemma follows.
\end{proof}

Now we are in a position to prove Theorem \ref{Theorem main} for the
remaining case.

\bigskip

{\sl Proof of Theorem \ref{Theorem main} for the case that $G=
Z_{pq}$ with $q\geq 2p+3.$}

Suppose $\ell_0\geq \lfloor2\sqrt{q-2}\rfloor$. By Lemma \ref{lemma
three facts} (iii), we have $\sum(A_0)=H=\langle A_0\rangle$. By
Observation \ref{Observation cap one subgroup}, we have
$A_0=H\setminus \{0\}$ and so $|\bigcup\limits_{i=1}^k
A_i|=|A|-|A_0|=p-2.$ It follows from Lemma \ref{Lemma fact on
sequence sums} that there exists an element $g\in G\setminus H$ such
that $\bigcup\limits_{i=1}^k A_i\subseteq (g+H)\cup (-g+H),$ we are
done. Therefore, we may assume that
\begin{equation}\label{equation |A0|<cr(zq)} \ell_0\leq
\lfloor2\sqrt{q-2}\rfloor-1,
\end{equation}
equivalently, \begin{equation}\label{equation sum l1+...+lk}
\sum\limits_{i=1}^k\ell_i\geq
(p+q-3)-(\lfloor2\sqrt{q-2}\rfloor-1)\geq p+3.
\end{equation}

\textbf{Claim A. } $\ell_i\leq p+1$ for all $i\in [1,k]$.

Assume to the contrary that $\ell_1\geq p+2$.  It follows from
\eqref{equation sum l1+...+lk} that (i) $\ell_1\geq p+3$ or (ii)
$\sum\limits_{i=2}^k \ell_i\geq 1$. If (i) holds, then
$|X_1^{\Delta}|=p$ and so $\alpha_1=p$. If (ii) holds, it is easy to
see  $|\sum\limits_{i=2}^{m_3} X_i+\sum\limits_{i=m_2+1}^k
X_i+Y|\geq 2$, by Lemma \ref{Lemma folk},
$\alpha_1=|X_1^{\Delta}+\sum\limits_{i=2}^{m_3}
X_i+\sum\limits_{i=m_2+1}^k X_i+Y|=p$. Since
$\ell_1+\ell_0+\epsilon(\ell_0)\geq p+2+2\geq 7$, applying Lemma
\ref{Lemma alpha and beta} with $n=1$, we derive a contradiction.
 This proves Claim A. \qed

\textbf{Claim B. } $\ell_0+\epsilon(\ell_0)-\delta(r_2)\leq 2$.

Assume to the contrary that $\ell_0+\epsilon(\ell_0)-\delta(r_2)\geq
3$. By Lemma \ref{Lemma alpha and beta}, we have $\alpha_0<p,$ and
thus, by \eqref{equation |A0|<cr(zq)} and Lemma \ref{Lemma alpha n},
\begin{align*}
p-1&\geq q-\ell_0-3 \\
&\geq q-(\lfloor2\sqrt{q-2}\rfloor-1)-3\\
&\geq \frac{q+1}{2}-3\\
&\geq \frac{(2p+3)+1}{2}-3 \\
&= p-1, \\
\end{align*}
which implies
\begin{equation}\label{equation q-..=q+1/2}
q-(\lfloor2\sqrt{q-2}\rfloor-1)=\frac{q+1}{2},
\end{equation} and
\begin{equation}\label{equation q=2p+3}
q=2p+3.
\end{equation}

By \eqref{equation q=2p+3}, we have $p\neq 3$ and so $q\geq 13$.
Thus, we check that
$$q-(\lfloor2\sqrt{q-2}\rfloor-1)>\frac{q+1}{2},$$ a contradiction
with \eqref{equation q-..=q+1/2}. This proves Claim B. \qed

By Claim B, we have $$\ell_0\leq 3.$$ Observe that
\begin{equation}\label{equation ell1+...ellk}
\sum\limits_{i=1}^k\ell_i=q+p-3-\ell_0\geq
(2p+3)+p-3-\ell_0=3(p-1)+(3-\ell_0).
\end{equation}

Suppose $\delta(r_2)=0$. By Claim B, we have $\ell_0\leq 2.$ By
\eqref{equation ell1+...ellk}, we have $\ell_1\geq 4$.  By Lemma
\ref{Lemma alpha n}, we check that $\alpha_1=p$. By Lemma \ref{Lemma
alpha and beta}, we derive a contradiction. Therefore,
\begin{equation}\label{equation r2>0}
\delta(r_2)=1.
\end{equation}

Suppose $\ell_0=3$. Since $\ell_k\leq 2$, it follows from
\eqref{equation ell1+...ellk} that  $\ell_1\geq 4$. Applying Lemma
\ref{Lemma alpha n} with $n=1$ and Lemma \ref{Lemma alpha and beta},
we derive a contradiction. Hence,
\begin{equation}\label{equation ell0 leq 2}
\ell_0\leq 2.
\end{equation}

Suppose $\ell_1\geq 5$. By Lemma \ref{Lemma alpha n}, we check that
$\alpha_1=p$. By Lemma \ref{Lemma alpha and beta}, we derive a
contradiction. Therefore,
\begin{equation}\label{equation ell 1 leq 4}
\ell_1\leq 4.
\end{equation}

By \eqref{equation ell1+...ellk}, \eqref{equation r2>0},
\eqref{equation ell0 leq 2} and \eqref{equation ell 1 leq 4}, we
conclude that $$\ell_1=\ell_2=4$$ and so $p\geq 5.$ It follows that
$$q=|A|-p+3\leq \ell_0+\sum\limits_{i=1}^3\ell_i+\ell_4-2\leq
2+3\times 4+2-2=14,$$ which implies $$q=13$$ and $$p=5.$$

Noting that $\ell_0\leq 2$, $\ell_1=\ell_2=4$, $\ell_4=2$ and
$\ell_3=5-\ell_0\geq 3$, we shall close this proof by deriving a
contradiction in the following.

Suppose that $\ell_3\geq 4$ or $\varphi(a_3)\neq-\varphi(a_4)$. Let
$X_4^{\blacklozenge}=\{0,\varphi(a_4),2\varphi(a_4)\}$. By Lemma
\ref{Lemma folk} and Lemma \ref{Lemma Vosper}, we have
 $|X_3+X_4^{\blacklozenge}|=p$, and so
 $|X_1^{\Delta}+X_2^{\Delta}+X_3+X_4^{\blacklozenge}|=p$. This implies that for
 every element $\bar{g}\in G\diagup H$, there exists a representation
 with coefficient $f_1=f_1(\bar{g}),\ldots,f_4=f_4(\bar{g})$,
 where $f_1=f_2=2$, $f_3\in [1,\ell_3-1]$ and $f_4\in
 [0,\ell_4]$. It follows from  Lemma \ref{Lemma Cauchy} and Lemma \ref{lemma three facts} (iii) that
\begin{align*}
&\ \ \ \ |\sum\nolimits^{\circ}(A_0)+\sum\limits_{i=1}^4
\sum\nolimits_{f_i}(\widetilde{A_i})|\\
&\geq|\sum\nolimits^{\circ}(A_0)+\sum\limits_{i=1}^3 \sum\nolimits_{f_i}(\widetilde{A_i})|\\
&\geq\min(q,|\sum\nolimits^{\circ}(A_0)|+\sum\limits_{i=1}^2|\sum\nolimits_{f_i}(\widetilde{A_i})|+|\sum\nolimits_{f_3}(\widetilde{A_3})|-3)\\
&\geq\min(q,(2\ell_0-1+\epsilon(\ell_0))+(\ell_1+1)+(\ell_2+1)+\ell_3-3)\\
&\geq\min(q,(2\ell_0-1+\epsilon(\ell_0))+5+5+(5-\ell_0)-3)\\
&=\min(q,\ell_0+\epsilon(\ell_0)+11)\\
&=q.
\end{align*}
This implies that $\sum(A)=G$, a contradiction. Therefore,
$$\ell_3=3$$ and $$\varphi(a_3)=-\varphi(a_4),$$ which implies
$$\varphi(a_2)\notin \{\varphi(a_3),-\varphi(a_3)\}.$$ By Lemma
\ref{Lemma Vosper}, we have
 $|X_2+X_3|\geq \min(p,|X_2|+|X_3|)=p$, and so
 $|X_1^{\Delta}+X_2+X_3+X_4|=p$. This implies that for
 every element $\bar{g}\in G\diagup H$, there exists a representation
 with coefficient $f_1=f_1(\bar{g}),\ldots,f_4=f_4(\bar{g})$,
 where $f_1=2$, $f_2\in [1,\ell_2-1]$, $f_3\in [1,\ell_3-1]$ and $f_4=1$. It follows
Lemma \ref{Lemma Cauchy} and Lemma \ref{lemma three facts} (iii)
that
\begin{align*}
&\ \ \ \ |\sum\nolimits^{\circ}(A_0)+\sum\limits_{i=1}^4
\sum\nolimits_{f_i}(\widetilde{A_i})|\\
&\geq\min(q,|\sum\nolimits^{\circ}(A_0)|+|\sum\nolimits_{f_1}(\widetilde{A_1})|+\sum\limits_{i=2}^4
|\sum\nolimits_{f_i}(\widetilde{A_i})|-4)\\
&\geq\min(q,(2\ell_0-1+\epsilon(\ell_0))+(\ell_1+1)+\ell_2+\ell_3+\ell_4-4)\\
&\geq\min(q,(2\ell_0-1+\epsilon(\ell_0))+5+4+(5-\ell_0)+2-4)\\
&\geq\min(q,\ell_0+\epsilon(\ell_0)+11)\\
&=q.
\end{align*}
This implies that $\sum(A)=G$, a contradiction.  \qed

\section{Conclusion}

We first give examples to show that the conclusion of Theorem
\ref{Theorem main} does not hold true for the group $G=Z_{pq}$ with
$q<2p+3.$

\begin{example}\label{example 1}\  Let $p,q$ be two odd prime numbers with
$p+\lfloor2\sqrt{p-2}\rfloor+1<q<2p+3$, and let $G=Z_{pq}$. Let $K$
be the subgroup of $G$ of order $p$. Let $A\subseteq K\cup (g+K)\cup
(-g+K)$ be a subset of $G$ with $|A|={\rm cr}(G)-1=p+q-3<3p$ and
$A\cap K=K\setminus \{0\}$, where $g\in G\setminus K$.
\end{example}

\begin{example}\label{example 2}\  Let $p,q$ be two odd prime numbers with
$p<q\leq p+\lfloor2\sqrt{p-2}\rfloor+1$, and let $G=Z_{pq}$. Let
$A=\{\pm g,\pm 2g,\ldots,\pm \frac{p+q-2}{2}g\}$ with $|A|={\rm
cr}(G)-1=p+q-2$, where $g\in G$ and ${\rm ord}(g)=pq$.
\end{example}

It is easy to check that $\sum(A)\neq G$ in both Example
\ref{example 1} and Example \ref{example 2}. However, we don't
observe any other counterexamples. Therefore, we conjecture the
following

\begin{conj}\ Let $p,q$ be two odd prime numbers with
$p+\lfloor2\sqrt{p-2}\rfloor+1<q<2p+3$, and let $A$ be a subset of
$G=Z_{pq}$ with $|A|={\rm cr}(G)-1=p+q-3$ and $0\notin A$. Then $A$
contains a complete subset.
\end{conj}

\begin{conj}\ Let $p,q$ be two odd prime numbers with
$p<q\leq p+\lfloor2\sqrt{p-2}\rfloor+1$, and let $A$ be a subset of
$G=Z_{pq}$ with $|A|={\rm cr}(G)-1=p+q-2$ and $0\notin A$. Then
there exists an element $g\in G$ of order $pq$ such that $A=\{\pm
g,\pm 2g,\ldots,\pm \frac{p+q-2}{2}g\}$
\end{conj}

We remark that, as shown in \cite{Vu2}, in most finite abelian
groups $G$, the comparatively large set $A$ fails to span the whole
group $G$ just because that most elements of $A$ concentrate in some
proper subgroup of $G$, i.e., $A$ contains a complete subset.
Therefore, we formulate the following theorem which is easy to
proved to be equivalent to Theorem B and Theorem \ref{Theorem main}
by Observation \ref{Observation cap one subgroup} and Lemma
\ref{Lemma fact on sequence sums}.

{\bf Theorem C.}\  {\sl Let $G$ be a finite abelian group, and let
$p$ be the smallest prime dividing $|G|$. Let $A$ be a subset of
$G\setminus \{0\}$ with $|A|={\rm cr}(G)-1$. Then $A$ contains a
complete subset provided that $G$ is one of the following types:

1. $|G|$ is an even number no less than $36$.

2. $\frac{|G|}{p}$ is composite and
$$\begin{array}{llll} & \frac{|G|}{p}\geq\left \{\begin{array}{llll}
               62, & \mbox{ if } \ \ p=3;\\
               7p+3, & \mbox{ if } \ \ p\geq 5.\\
              \end{array}
           \right. \\
\end{array}$$

3. $|G|$ is a product of two odd prime numbers $p,q$ with $q\geq
2p+3$.

\bigskip

\noindent {\bf Acknowledgement.} The authors are grateful to the the
referee for helpful suggestions and comments.

\end{document}